\documentclass[a4paper,12pt]{article}
 
 \usepackage{latexsym,amssymb,amsmath}  
 
 \newtheorem{rem}{Remark}

\newcommand{\uz}{\underline z}
\newcommand{\uzd}{\underline z^{\dagger}}
\newcommand{\puz}{\partial_{\underline z}}
\newcommand{\puzd}{\partial_{\underline z^{\dagger}}}

\begin{document}

\title{Special functions and systems in Hermitian Clifford analysis}

\author{Nele De Schepper$^{1}$, Dixan Pe\~na Pe\~na$^{2}$ and Frank Sommen$^{3}$}

\date{\normalsize{Clifford Research Group, Department of Mathematical Analysis\\Ghent University\\Galglaan 2, 9000 Gent, Belgium}\\\vspace{0.1cm}
\small{$^{1}$e-mail: nds@cage.ugent.be\\
$^{2}$e-mail: dpp@cage.ugent.be\\
$^{3}$e-mail: fs@cage.ugent.be}}

\maketitle

\begin{abstract}
\noindent 
In this paper we study some new special functions that arise naturally within the framework of Hermitian Clifford analysis, which concerns the study of Dirac-like systems in several complex variables. In particular we focus on Hermite polynomials, Bessel functions and generalized powers. We also derive a Vekua system for solutions of Hermitian systems in axially symmetric domains.\vspace{0.2cm}\\
\textit{Keywords}: Hermitian Clifford analysis; Cauchy-Kowalevski extension; Bessel functions; Laguerre polynomials.\vspace{0.1cm}\\
\textit{Mathematics Subject Classification}: 30G35, 33C10, 33C45.
\end{abstract}

%%%%%%%%%%%%%%%%%%%%%%%%%%%%%%%%%%%%%%%%%%%%%%%%%%%%%%%%%%%%%%%%%%%%%%%%%%%%%%%%%%%%%%%%%%%%%%%%%%%%%%%%%%%%%%%%%%%%%%%%%%%%%%%%
\section{Introduction}
Clifford analysis deals with partial differential operators that arise naturally within the context of a Clifford algebra. The operator that is studied most is the generalized Cauchy-Riemann operator $\partial_{X_0}+\partial_{\underline{X}}$, $\partial_{\underline{X}}=\sum_{j=1}^me_j \partial_{X_j}$ being the Dirac operator and $e_1,\dots,e_m$ the generators of the Clifford algebra satisfying the defining relations $e_je_k+e_ke_j=-2\delta_{jk}$. Solutions of $(\partial_{X_0}+\partial_{\underline{X}})f=0$ are called monogenic functions in $\mathbb R^{m+1}$. One of the main basic properties is the Cauchy-Kowalevski extension (CK-extension) theorem which says that a given real analytic function $g(\underline X)$, $\underline X\in\mathbb R^m$, admits a unique monogenic extension $f(X_0,\underline X)$ is some domain of $\mathbb R^{m+1}$ that is given by
\begin{equation*}
f(X_0,\underline X)=\sum_{k=0}^\infty\frac{(-X_0)^k}{k!}\,\partial_{\underline X}^kg(\underline X).
\end{equation*}
In other words, the restriction operator $f\rightarrow f\vert_{X_0=0}$ is injective and also surjective. This CK-extension leads to the construction of numerous special monogenic functions depending on the choice of the initial function $g(\underline X)$ and it gives rise to special classes of polynomials and other special functions on the way. In particular in \cite{S2} was constructed Clifford-Bessel functions, monogenic generalized powers and Clifford Hermite polynomials arising from the CK-extension of the Gaussian distribution (the monogenic extension of the Gaussian distribution for $m$ odd has been obtained in closed form in \cite{DPSo} using Fueter's theorem). The above generalized Cauchy-Riemann system is invariant under the orthogonal group (the spin group $\text{Spin}(m)$ to be more precise).

There also exists a so-called Hermitian monogenic system that can be defined in several complex variables (see e.g. \cite{H6A, H6A2, BHS, R-CSS,SS}). It is invariant under the action of the unitary group and it contains the holomorphic functions in several complex variables as a subclass. To define that system we consider several complex variables like $z_0,z_1,\dots,z_n$ and their conjugates $\overline z_0,\overline z_1,\dots,\overline z_n$ and study solutions of the system
\[(f_0^\dagger\partial_{z_0}+f_1^\dagger\partial_{z_1}+\dots+f_n^\dagger\partial_{z_n})f(z_j,\overline z_j)=(f_0\partial_{\overline z_0}+f_1\partial_{\overline z_1}+\dots+f_n\partial_{\overline z_n})f(z_j,\overline z_j)=0,\]
called Hermitian monogenic ($h$-monogenic) functions, whereby the Clifford generators $f_j,f_j^\dagger$ satisfy the relations for a Witt basis
\[f_jf_k=-f_kf_j,\quad f_j^\dagger f_k^\dagger=-f_k^\dagger f_j^\dagger,\quad f_jf_k^\dagger+f_k^\dagger f_j=\delta_{jk}.\]
Like in the orthogonal setting one may wonder about a CK-extension theorem for $h$-monogenic functions; it has been studied in \cite{BHLS}. However, this time the restriction operator $f\rightarrow f\vert_{z_0=0}$ is no longer injective and it is also not surjective so that e.g. the Gaussian $\exp(-\frac{1}{2}\sum_{j=1}^nz_j\overline z_j)$ admits no CK-extension. Yet in our paper \cite{HCHP-AACA} we were able to construct an analogue of Hermite polynomials for the unitary group ($h$-monogenic setting). This leads to the consideration of a weaker form of the $h$-monogenic system, containing the $h$-monogenic functions as a subclass, which is just big enough for the restriction operator $f\rightarrow f\vert_{z_0=0}$ to be surjective.

The main topic in this paper is to study the CK-extension theorem for this sub-system of the $h$-monogenic system, showing that its solutions (called \textit{$h$-submonogenic functions}) are determined by their Cauchy data and to solve the system explicitly for the Gaussian and for other special functions as well. In this way we obtain Hermite polynomials, generalized powers and Bessel functions for the  $h$-submonogenic system. We also derive a Vekua-type system (see e.g. \cite{Krav,Ve}) that describes all axially symmetric solutions to this system. It is the Hermitian analogue of the Vekua system for the so-called axial monogenic functions (see \cite{LB,S1,S2}). 

In Sections \ref{sect2} and \ref{sect3} we recall the details concerning the Hermitian monogenic system and introduce the sub-system that will lead to a surjective CK-extension theorem. In Section \ref{sect4} we also discuss the constraints that have to be imposed on the Cauchy data for the CK-extension to be $h$-monogenic. Then we derive the Hermite polynomials from the CK-extension of the Gaussian. In the next section we establish the Vekua system for axially symmetric solutions and apply this to the construction of generalized powers. In the last section we study solutions of exponential type, leading to Hermitian Bessel functions, which have applications in Fourier analysis.
%%%%%%%%%%%%%%%%%%%%%%%%%%%%%%%%%%%%%%%%%%%%%%%%%%%%%%%%%%%%%%%%%%%%%%%%%%%%%%%%%%%%%%%%%%%%%%%%%%%%%%%%%%%%%%%%%%%%%%%%%%%%%%%%%%%
\section{Hermitian Clifford Analysis}\label{sect2}
In $m$-dimensional Euclidean space, Clifford analysis (see e.g. \cite{BDS, DSS, GM, GS}) focusses on the null solutions of various partial differential operators arising within the Clifford algebra language, the most important one being the Dirac operator $\partial_{\underline{X}}=\sum_{j=1}^me_j \partial_{X_j}$, which is the Fischer dual of the Clifford vector variable $\underline{X}=\sum_{j=1}^mX_je_j$, and the null solutions of which are called monogenic functions. Here $(e_1,\ldots,e_m)$ forms the usual orthonormal basis for the real vector space $\mathbb{R}^m$, equipped with a bilinear form of signature $(0,m)$ and underlying the construction of the real Clifford algebra $\mathbb{R}_{0,m} = \oplus_{k=0}^m \mathbb{R}_{0,m}^k$ where $\mathbb{R}_{0,m}^k$ denotes the subspace of $k$-vectors, spanned by the products of $k$ different basis vectors. This Clifford or geometric product is governed by the well-known noncommutative rules $e_j^2=-1$, $e_j e_k + e_k e_j=0$, $j \not= k=1,\ldots , m$. We refer to this setting as the orthogonal case, since the fundamental group leaving the Dirac operator $\partial_{\underline{X}}$ invariant is the special orthogonal group $\text{SO}(m)$, which is doubly covered by the $\text{Spin}(m)$ group of the Clifford algebra $\mathbb{R}_{0,m}$.

When allowing for complex scalars, the same set of generators as above $(e_1, \ldots,e_m)$, still satisfying the defining relations $e_j e_k + e_k e_j = -2 \delta_{jk}$, $j,k=1,\ldots,m$, may in fact also generate the complex Clifford algebra $\mathbb{C}_m$, where moreover we take the dimension to be even, say $m=2n$, for intrinsic reasons, see e.g. \cite{H6A,H6A2}. As $\mathbb{C}_{2n}$ is the complexification of the real Clifford algebra $\mathbb{R}_{0,2n}$, i.e. $\mathbb{C}_{2n} = \mathbb{R}_{0,2n} \oplus i \mathbb{R}_{0,2n}$, any complex Clifford number $\lambda \in \mathbb{C}_{2n}$ may be written as $\lambda = a+ib$, $a,b \in \mathbb{R}_{0,2n}$, leading to the definition of the Hermitian conjugation:
\begin{displaymath}
\lambda^{\dagger} = (a+ib)^{\dagger} = \overline{a} -i \overline{b},
\end{displaymath}
where the bar denotes the usual conjugation in $\mathbb{R}_{0,2n}$, i.e. the unique anti-involution for which $\overline{e}_j = -e_j$, $j=1,\ldots,2n$. This Hermitian conjugation then gives rise to a Hermitian inner product and its associated norm on $\mathbb{C}_{2n}$ given by
\begin{displaymath}
(\lambda, \mu) = \lbrack \lambda^{\dagger} \mu \rbrack_0, \qquad \vert \lambda \vert = \sqrt{\lbrack \lambda^{\dagger} \lambda \rbrack_0},
\end{displaymath}
where $\lbrack \ \rbrack_0$ denotes the scalar part.

The above framework will be referred to as the Hermitian Clifford setting, as opposed to the traditional orthogonal Clifford setting. For the complex Clifford algebra $\mathbb{C}_{2n}$ we consider the so-called Witt basis
\begin{displaymath}
f_j = \frac{1}{2} (e_j-ie_{n+j}), \qquad f_j^{\dagger} = - \frac{1}{2} (e_j + i e_{n+j}), \ \ j=1,\ldots,n.
\end{displaymath}
These Witt basis elements are isotropic
\begin{displaymath}
f_j^2=(f_j^{\dagger})^2 = 0, \quad j=1,\ldots,n
\end{displaymath}
and satisfy the Grassmann identities
\begin{displaymath}
f_j f_k + f_k f_j = f_j^{\dagger} f_k^{\dagger} + f_k^{\dagger} f_j^{\dagger}=0, \quad j,k=1,\ldots, n
\end{displaymath}
as well as the duality identities
\begin{displaymath}
f_j f_k^{\dagger} + f_k^{\dagger} f_j=\delta_{jk}, \quad j,k=1, \ldots, n.
\end{displaymath}
The Clifford algebra $\mathbb{C}_{2n}$ is isomorphic to the full matrix algebra $\mathbb C(2^{2n})=\text{End}(S)$, whereby the so-called spinor space $S$ is the $2^n$-dimensional minimal left ideal $S=\mathbb{C}_{2n}I$, $I$ being the primitive idempotent
\[I=f_1f_1^\dagger\dots f_nf_n^\dagger.\] 
As clearly $f_jI=0$, the spinor space is given by $S=\Lambda_nI$, $\Lambda_n$ being the Grassmann algebra generated by the elements $f_1^\dagger,\dots,f_n^\dagger$. Grassmann algebras may be seen as fermionic polynomial algebras generated by anticommuting fermionic variables such as $f_1^\dagger,\dots,f_n^\dagger$. The corresponding fermionic derivatives $\partial_{f_1^\dagger},\dots,\partial_{f_n^\dagger}$ are clearly the elements $f_1,\dots,f_n$ and the defining relations for $\mathbb{C}_{2n}$ in terms of the Witt basis are the fermionic equivalent of the Weyl relations. The fermionic Euler operator $\sum_{j=1}^nf_j^\dagger\partial_{f_j^\dagger}$ is the Clifford number $\beta=\sum_{j=1}^nf_j^\dagger f_j$ and its eigenspaces are the spaces $\Lambda_n^{\ell}I$ whereby $\Lambda_n^{\ell}$ is the space of fermionic polynomials of degree $\ell$ spanned by the $\ell$-vectors $f_A^\dagger=f_{\alpha_1}^\dagger\dots f_{\alpha_{\ell}}^\dagger$. The corresponding eigenvalues are given by $\beta\Lambda_n^{\ell}I=\ell\Lambda_n^{\ell}I$ so that in fact $\beta$ satisfies the characteristic equation 
\[\beta(\beta-1)\dots(\beta-n)=0.\]
Using this Witt basis, the vector $(X_1, \ldots, X_{2n}) = (x_1, \ldots, x_n, y_1, \ldots , y_n)$ in $\mathbb{R}^{2n}$ is identified with the Clifford vector
\begin{displaymath}
\underline{X} = \sum_{j=1}^n (x_j e_j + y_j e_{n+j}) = \sum_{j=1}^n f_j z_j - \sum_{j=1}^n f_j^{\dagger} \overline{z}_j,
\end{displaymath}
where the complex variables $z_j=x_j + i y_j$ and their complex conjugates $ \overline{z}_j= x_j - i y_j$, $j=1,\ldots,n$ have been introduced. Defining the Hermitian vector variable $\underline{z}$ and its Hermitian conjugate $\underline{z}^{\dagger}$ by 
\begin{displaymath}
\uz = \sum_{j=1}^n f_j z_j, \qquad \qquad \uzd = \sum_{j=1}^n f_j^{\dagger} \overline{z}_j,
\end{displaymath}
the Clifford vector $\underline{X}$ clearly decomposes as $\underline{X} = \uz - \uzd$. This also gives rise to the decomposition of the traditional Dirac operator 
\begin{displaymath}
\partial_{\underline{X}} = \sum_{j=1}^n (e_j\partial_{x_j} +e_{n+j}\partial_{y_j}) = 2 \sum_{j=1}^n (f_j \partial_{\overline{z}_j} - f_j^{\dagger} \partial_{z_j}) = 2 (\puzd- \puz)
\end{displaymath} 
in terms of the Hermitian Dirac operators
\begin{displaymath}
\puz = \sum_{j=1}^n  f_j^{\dagger} \partial_{z_j}, \qquad \qquad \puzd = (\puz)^{\dagger} = \sum_{j=1}^n  f_j \partial_{\overline{z}_j}
\end{displaymath}
involving the classical Cauchy-Riemann operators and their complex conjugates in the complex $z_j$ planes, i.e. $\partial_{z_j} = \frac{1}{2} (\partial_{x_j} - i \partial_{y_j})$ and $\partial_{\overline{z}_j} = \frac{1}{2} (\partial_{x_j} + i \partial_{y_j})$, $j=1,\ldots,n$.\\
On account of the isotropy of the Witt basis elements, the Hermitian vector variables and Dirac operators are isotropic as well, from which it directly follows that the Laplacian $\Delta_{2n} = - \partial_{\underline{X}}^2$ in $\mathbb{R}^{2n}$ allows for the decomposition $\Delta_{2n} = 4 (\puz \puzd + \puzd \puz)$, which is the dual expression of
\begin{displaymath}
\vert \uz \vert^2 = \vert \uzd \vert^2 = \uz\,\uzd + \uzd \uz.
\end{displaymath}
Hermitian Clifford analysis (see e.g. \cite{H6A, H6A2, BHS, R-CSS,SS}) then focusses on the null-solutions of both Hermitian Dirac operators obtained. Indeed, a continuously differentiable function $g$ on $\mathbb{R}^{2n}$ with values in $\mathbb{C}_{2n}$ is called a $h$-monogenic function if and only if it satisfies the system
\begin{equation}\label{hsystem}
\puz g = 0 = \puzd g,
\end{equation}
or in other words, iff it is a simultaneous null solution of the Hermitian Dirac operators $\puz$ and $\puzd$.
%%%%%%%%%%%%%%%%%%%%%%%%%%%%%%%%%%%%%%%%%%%%%%%%%%%%%%%%%%%%%%%%%%%%%%%%%%%%%%%%%%%%%%%%%%%%%%%%%%%%%%%%%%%%%%%%%%%%%%%%%%%%%%%%
\section{Introducing our system}\label{sect3}
Let us consider the $h$-monogenic system (\ref{hsystem}) in dimension $m=2n+2$
\begin{equation}\label{hms}
\left\{\begin{aligned}
(f_0^{\dagger}\partial_{z_0}+\partial_{\underline z})f&=0\\
(f_0\partial_{\overline{z}_0}+\partial_{\underline z^{\dagger}})f&=0
\end{aligned}\right.
\end{equation}
where $f:\Omega\subset\mathbb R^{2n+2}\rightarrow\mathbb{C}_{2n+2}$ is a continuously differentiable function. If we multiply the first equation by $f_0^{\dagger}f_0$ (resp. $f_0f_0^{\dagger}$), we get
\[(f_0^{\dagger}\partial_{z_0}+f_0^{\dagger}f_0\partial_{\underline z})f=0\;(\text{resp.}\;f_0f_0^{\dagger}\partial_{\underline z}f=0).\] 
These last two equations are clearly equivalent to the first equation of (\ref{hms}). In a similar way, we can also show that the second equation of (\ref{hms}) is equivalent to 
\[(f_0\partial_{\overline z_0}+f_0f_0^{\dagger}\partial_{\underline z^{\dagger}})f=f_0^{\dagger}f_0\partial_{\underline z^{\dagger}}f=0.\] 
In other words, $f$ is a solution of the $h$-monogenic system (\ref{hms}) if and only if  
\begin{equation}\label{hmsF}
\left\{\begin{aligned}
(f_0^{\dagger}\partial_{z_0}+f_0^{\dagger}f_0\partial_{\underline z})f&=0\\
f_0f_0^{\dagger}\partial_{\underline z}f&=0\\
(f_0\partial_{\overline z_0}+f_0f_0^{\dagger}\partial_{\underline z^{\dagger}})f&=0\\
f_0^{\dagger}f_0\partial_{\underline z^{\dagger}}f&=0.
\end{aligned}\right.
\end{equation} 
In this paper we shall focus on the following system of equations
\begin{equation}\label{whms}
\left\{\begin{aligned}
(f_0^{\dagger}\partial_{z_0}+f_0^{\dagger}f_0\partial_{\underline z})f&=0\\
(f_0\partial_{\overline z_0}+f_0f_0^{\dagger}\partial_{\underline z^{\dagger}})f&=0
\end{aligned}\right.
\end{equation} 
which is a sub-system of the $h$-monogenic system (\ref{hmsF}), solutions of which will hence be called $h$-submonogenic functions. 

While the $h$-monogenic system in $\mathbb C^{n+1}$ is invariant under the action of the unitary group $\text{U}(n+1)$, the $h$-submonogenic system (\ref{whms}) no longer enjoys this invariance and it depends on the choice of the special complex direction $z_0$. But the system is still invariant under the action of the unitary subgroup $\text{U}(n)$ of $\text{U}(n+1)$. Moreover, clearly the $h$-submonogenic system (\ref{whms}) is equivalent to the inhomogeneous $h$-monogenic system
\begin{equation}\label{inhhms}
\left\{\begin{aligned}
(f_0^{\dagger}\partial_{z_0}+\partial_{\underline z})f&=f_0g\\
(f_0\partial_{\overline{z}_0}+\partial_{\underline z^{\dagger}})f&=f_0^{\dagger}h
\end{aligned}\right.
\end{equation}
 whereby $g=f_0^{\dagger}\partial_{\underline z}f$, $h=f_0\partial_{\underline z^{\dagger}}f$. This puts things back into the $\text{U}(n+1)$-invariant setting of the inhomogeneous $h$-monogenic system, whereby the symmetry breaking comes from the specific choice of the right-hand sides ($f_0g$, $f_0^{\dagger}h$) of (\ref{inhhms}).

We will now study the Cauchy-Kowalevski extension problem for the $h$-submonogenic system. That is, given $g:\underline\Omega\subset\mathbb R^{2n}\rightarrow\mathbb{C}_{2n+2}$ , find a solution $f$ of (\ref{whms}) such that $f\vert_{z_0=0}=g$.

First note that any $\mathbb{C}_{2n+2}$-valued function $f$ may be uniquely written into the form
\[f=A+f_0B+f_0^{\dagger}C+f_0^{\dagger}f_0D,\]
where $A$, $B$, $C$, $D$ are $\mathbb{C}_{2n}$-valued functions.

A direct computation then yields
\begin{align*}
(f_0^{\dagger}\partial_{z_0}+f_0^{\dagger}f_0\partial_{\underline z})f&=f_0^{\dagger}(\partial_{z_0}A-\partial_{\underline z}C)+f_0^{\dagger}f_0\big(\partial_{z_0}B+\partial_{\underline z}(A+D)\big)\\
(f_0\partial_{\overline z_0}+f_0f_0^{\dagger}\partial_{\underline z^{\dagger}})f&=f_0\big(\partial_{\overline z_0}(A+D)-\partial_{\underline z^{\dagger}}B\big)+f_0f_0^{\dagger}(\partial_{\underline z^{\dagger}}A+\partial_{\overline z_0}C).
\end{align*}
Therefore the $h$-submonogenic system (\ref{whms}) is equivalent to the following systems of equations
\begin{equation}\label{whms1}
\left\{\begin{aligned}
\partial_{z_0}A&=\partial_{\underline z}C\\
\partial_{\underline z^{\dagger}}A&=-\partial_{\overline z_0}C
\end{aligned}\right.
\end{equation}
\begin{equation}\label{whms2}
\left\{\begin{aligned}
\partial_{z_0}B&=-\partial_{\underline z}(A+D)\\
\partial_{\underline z^{\dagger}}B&=\partial_{\overline z_0}(A+D).
\end{aligned}\right.
\end{equation}
\begin{rem}
Note that if moreover, we want $f$ to be a $h$-monogenic function, then the following extra conditions should be satisfied
\[\partial_{\underline z}A=\partial_{\underline z^{\dagger}}(A+D)=\partial_{\underline z}B=\partial_{\underline z^{\dagger}}C=0.\] 
\end{rem}
%%%%%%%%%%%%%%%%%%%%%%%%%%%%%%%%%%%%%%%%%%%%%%%%%%%%%%%%%%%%%%%%%%%%%%%%%%%%%%%%%%%%%%%%%%%%%%%%%%%%%%%%%%%%%%%%%%%%%%%%%%%%%%%%%%%%%%%%%%%%
\section{Cauchy-Kowalevski extension problem}\label{sect4}
\subsection{Power series method}\label{PSMet}
Assume
\begin{equation}\label{middle}
\begin{split}
A&=\sum_{k=0}^\infty(z_0\overline z_0)^kA_k,\qquad\quad\qquad B=\sum_{k=0}^\infty z_0(z_0\overline z_0)^kB_k,\\
C&=\sum_{k=0}^\infty\overline z_0(z_0\overline z_0)^kC_k,\qquad\qquad D=\sum_{k=0}^\infty(z_0\overline z_0)^kD_k,
\end{split}
\end{equation}
where $A_k$, $B_k$, $C_k$, $D_k$ are $\mathbb{C}_{2n}$-valued continuously differentiable functions defined on $\mathbb R^{2n}$.
\begin{rem}
Note that $f\vert_{z_0=0} = A_0 + f_0^{\dagger} f_0 D_0$. This restriction uniquely determines $f$. Hence, the initial conditions for the solution (\ref{middle}) of the $h$-submonogenic system (\ref{whms}) are the functions $A_0$ and $D_0$.
\end{rem}

Substituting the expressions (\ref{middle}) in the systems (\ref{whms1}) and (\ref{whms2}), we get the recurrence relations
\begin{align}
A_{k+1}&=\frac{1}{k+1}\partial_{\underline z}C_k,&\label{rr1}
C_k&=-\frac{1}{k+1}\partial_{\underline z^{\dagger}}A_k,&k\ge0,\\\label{rr2}
B_k&=-\frac{1}{k+1}\partial_{\underline z}(A_k+D_k),&
A_{k+1}+D_{k+1}&=\frac{1}{k+1}\partial_{\underline z^{\dagger}}B_k,&k\ge0.
\end{align}
Then, from (\ref{rr1}) we obtain the solution
\begin{align*}
A_k&=\frac{(-1)^k}{(k!)^2}(\partial_{\underline z}\partial_{\underline z^{\dagger}})^kA_0,&k\ge1,\\
C_k&=\frac{(-1)^{k+1}}{(k+1)(k!)^2}\partial_{\underline z^{\dagger}}(\partial_{\underline z}\partial_{\underline z^{\dagger}})^kA_0,&k\ge0.
\end{align*}
Similarly, from (\ref{rr2}) we obtain that 
\begin{align*}
B_k&=\frac{(-1)^{k+1}}{(k+1)(k!)^2}\partial_{\underline z}(\partial_{\underline z^{\dagger}}\partial_{\underline z})^k(A_0+D_0),&k\ge0,\\
D_k&=\frac{(-1)^k}{(k!)^2}\Big[(\partial_{\underline z^{\dagger}}\partial_{\underline z})^kD_0+\Big((\partial_{\underline z^{\dagger}}\partial_{\underline z})^k-(\partial_{\underline z}\partial_{\underline z^{\dagger}})^k\Big)A_0\Big],&k\ge1.
\end{align*}
\begin{rem}
If we want the $h$-submonogenic CK-extension $f$ of $f\vert_{z_0=0} = A_0 + f_0^{\dagger} f_0 D_0$ to be $h$-monogenic, the initial conditions $A_0$ and $D_0$ should satisfy the extra constraints
\begin{displaymath}
\puz A_0 = \puzd (A_0 + D_0)=0.
\end{displaymath}
\end{rem}
%%%%%%%%%%%%%%%%%%%%%%%%%%%%%%%%%%%%%%%%%%%%%%%%%%%%%%%%%%%%%%%%%%%%%%%
\subsection{Clifford-Hermite polynomials}
Let us consider for example the initial condition $f\vert_{z_0=0} = e^{-\frac{\vert \uz \vert^2}{2}}$, i.e. $A_0 = e^{-\frac{\vert \uz \vert^2}{2}}$ and $D_0 = 0$. Note that $\puz A_0 \not=0$, hence $f=A+f_0B + f_0^{\dagger} C + f_0^{\dagger} f_0 D$ cannot be $h$-monogenic.

Expressing $A_k$, $B_k$, $C_k$ and $D_k$ in terms of powers of the Laplacian acting on the Gauss function, we obtain:
\[A_k=\frac{(-1)^k 4^{1-k}}{(k!)^2} (\puz \puzd) \Delta_{2n}^{k-1} e^{-\frac{\vert \uz \vert^2}{2}},\quad B_k= \frac{(-1)^{k+1} 4^{-k}}{(k+1)(k!)^2} \puz  \Delta_{2n}^{k} e^{-\frac{\vert \uz \vert^2}{2}},\]
\[C_k=\frac{(-1)^{k+1} 4^{-k}}{(k+1)(k!)^2} \puzd  \Delta_{2n}^{k} e^{-\frac{\vert \uz \vert^2}{2}},\]
\[D_k=\frac{(-1)^k 4^{1-k}}{(k!)^2} \left( \puzd \puz \Delta_{2n}^{k-1} e^{-\frac{\vert \uz \vert^2}{2}}- \puz \puzd \Delta_{2n}^{k-1} e^{-\frac{\vert \uz \vert^2}{2}} \right).\]
The above functions can be rewritten in terms of the \textit{Hermitian Clifford-Hermite polynomials} introduced in \cite{HCHP-AACA} by means of a Rodrigues type formula involving both Hermitian Dirac operators:
\begin{align*}
H_{2p+1}^{(1)}(\uz,\uzd) e^{-\frac{\vert \uz \vert^2}{2}} &= \puzd \Delta_{2n}^{p} e^{-\frac{\vert \uz \vert^2}{2}}, \qquad H_{2p+1}^{(2)}(\uz,\uzd) e^{-\frac{\vert \uz \vert^2}{2}} = \puz \Delta_{2n}^{p} e^{-\frac{\vert \uz \vert^2}{2}}\\
H_{2p+2}^{(3)}(\uz,\uzd) e^{-\frac{\vert \uz \vert^2}{2}} &= \puz \puzd \Delta_{2n}^{p} e^{-\frac{\vert \uz \vert^2}{2}}, \quad H_{2p+2}^{(4)}(\uz,\uzd) e^{-\frac{\vert \uz \vert^2}{2}} = \puzd \puz \Delta_{2n}^{p} e^{-\frac{\vert \uz \vert^2}{2}},
\end{align*}
$p\in\mathbb N_0$. These polynomials  can be written in terms of the Laguerre polynomials on the real line:
\[H_{2p+1}^{(1)}(\uz,\uzd)= (-1)^{p-1} 2^{p-1} p! \uz L_p^n \left( \frac{\vert \uz \vert^2}{2} \right)\]
\[H_{2p+1}^{(2)}(\uz,\uzd)= (-1)^{p-1} 2^{p-1} p! \uzd L_p^n \left( \frac{\vert \uz \vert^2}{2} \right)\]
\[H_{2p+2}^{(3)}(\uz,\uzd)= (-1)^{p-1} 2^{p-1} p! \left( \beta  L_p^n \left( \frac{\vert \uz \vert^2}{2} \right) - \frac{1}{2} \uzd \uz L_p^{n+1} \left( \frac{\vert \uz \vert^2}{2} \right) \right)\]
\[H_{2p+2}^{(4)}(\uz,\uzd)= (-1)^{p-1} 2^{p-1} p! \left( (n-\beta)  L_p^n \left( \frac{\vert \uz \vert^2}{2} \right) - \frac{1}{2} \uz\,\uzd L_p^{n+1} \left( \frac{\vert \uz \vert^2}{2} \right) \right)\]
whereby $\beta=\sum_{j=1}^n f_j^{\dagger} f_j$ is the fermionic Euler operator. Hence, the functions $A_k$, $B_k$, $C_k$ and $D_k$ take the following form:
\begin{align*}
A_k &= \frac{(-1)^k 4^{1-k}}{(k!)^2} e^{-\frac{\vert \uz \vert^2}{2}} H_{2k}^{(3)}(\uz,\uzd)\\
&= \frac{2^{-k}}{k k!} e^{-\frac{\vert \uz \vert^2}{2}} \left( \beta L_{k-1}^n \left( \frac{\vert \uz \vert^2}{2} \right) - \frac{1}{2} \uzd \uz L_{k-1}^{n+1} \left( \frac{\vert \uz \vert^2}{2} \right) \right), \ k \geq 1
\end{align*}
\[B_k= \frac{(-1)^{k+1} 4^{-k}}{(k+1) (k!)^2} e^{-\frac{\vert \uz \vert^2}{2}} H_{2k+1}^{(2)}(\uz,\uzd) = \frac{2^{-k-1}}{(k+1)!}  e^{-\frac{\vert \uz \vert^2}{2}} \uzd L_{k}^n \left( \frac{\vert \uz \vert^2}{2} \right), \ k \geq 0\]
\[C_k= \frac{(-1)^{k+1} 4^{-k}}{(k+1) (k!)^2} e^{-\frac{\vert \uz \vert^2}{2}} H_{2k+1}^{(1)}(\uz,\uzd) = \frac{2^{-k-1}}{(k+1)!}  e^{-\frac{\vert \uz \vert^2}{2}} \uz L_{k}^n \left( \frac{\vert \uz \vert^2}{2} \right), \ k \geq 0\]
\begin{align*}
D_k &= \frac{(-1)^k 4^{1-k}}{(k!)^2} e^{-\frac{\vert \uz \vert^2}{2}}  \left( H_{2k}^{(4)}(\uz,\uzd) - H_{2k}^{(3)}(\uz,\uzd) \right), \ k \geq 1\\
&= \frac{2^{-k}}{k k!} e^{-\frac{\vert \uz \vert^2}{2}} \left( k  L_{k}^n \left( \frac{\vert \uz \vert^2}{2} \right) -  (2\beta +k)  L_{k-1}^n \left( \frac{\vert \uz \vert^2}{2} \right) + \uzd \uz L_{k-1}^{n+1} \left( \frac{\vert \uz \vert^2}{2} \right) \right),  
\end{align*}
where in the last line we have used the relation
\begin{displaymath}
x L_{n-1}^{(\alpha+1)}(x) = (n+\alpha) L_{n-1}^{(\alpha)}(x) - n L_{n}^{(\alpha)}(x).
\end{displaymath}
%%%%%%%%%%%%%%%%%%%%%%%%%%%%%%%%%%%%%%%%%%%%%%%%%%%%%%%%%%%%%%%%%%%%%%%%%%%%%%%%%%%%%%%%%%%%%%%%%%%%%%%%%%%%%%%%%%%%%%%%%%%%%%%%%%%%%%%
\subsection{Double power series method}
We now generalize the power series method of the previous subsection by putting
\begin{align*}
A&=\sum_{k, \ell=0}^\infty z_0^k\,\overline z_0^{\ell} A_{k,\ell},&
B&=\sum_{k, \ell=0}^\infty z_0^k\,\overline z_0^{\ell} B_{k,\ell},\\
C&=\sum_{k, \ell=0}^\infty z_0^k\,\overline z_0^{\ell} C_{k,\ell},&
D&=\sum_{k, \ell=0}^\infty z_0^k\,\overline z_0^{\ell} D_{k,\ell},
\end{align*}
with $A_{k,\ell}$, $B_{k,\ell}$, $C_{k,\ell}$ and $D_{k,\ell}$ $\mathbb{C}_{2n}$-valued continuously differentiable functions defined on $\mathbb{R}^{2n}$. Note that $f\vert_{z_0=0} = A_{0,0} + f_0 B_{0,0} + f_0^{\dagger} C_{0,0} + f_0^{\dagger} f_0 D_{0,0}$.

Substitution of these expressions in the systems (\ref{whms1}) and (\ref{whms2}) now leads to the recurrence relations
\begin{equation}\label{recI}
A_{k+1,\ell} = \frac{1}{k+1} \puz C_{k, \ell}, \quad C_{k,\ell+1} = - \frac{1}{\ell+1} \puzd A_{k,\ell}, 
\end{equation}
and
\begin{equation}\label{recII}
B_{k+1,\ell} = - \frac{1}{k+1} \puz (A_{k,\ell} + D_{k,\ell}),\quad A_{k,\ell+1} + D_{k,\ell+1} = \frac{1}{\ell+1} \puzd B_{k,\ell}.
\end{equation}
Equations (\ref{recI}) yield
\begin{align*}
A_{k,0} &= \frac{1}{k} \puz C_{k-1,0}, \ \ k \geq 1\\
A_{k,\ell} &= (-1)^k \frac{(\ell-k)!}{k! \ell!} (\puz \puzd)^k A_{0,\ell-k}, \ \ 1 \leq k \leq \ell\\
A_{k, \ell} &= (-1)^{\ell} \frac{(k-1-\ell)!}{k! \ell!} \puz (\puzd \puz)^{\ell} C_{k-1-\ell,0}, \ \ 1 \leq \ell < k
\end{align*}
and 
\begin{align*}
C_{0,\ell} &= - \frac{1}{\ell} \puzd A_{0,\ell-1}, \ \ \ell \geq 1\\
C_{k,\ell} &= (-1)^{k+1} \frac{(\ell-1-k)!}{k! \ell!} \puzd (\puz \puzd)^k A_{0,\ell-1-k}, \ \ 1 \leq k < \ell\\
C_{k, \ell} &= (-1)^{\ell} \frac{(k-\ell)!}{k! \ell!} (\puzd \puz)^{\ell} C_{k-\ell,0}, \ \ 1 \leq \ell \leq k,
\end{align*}
while from (\ref{recII}) we obtain
\begin{align*}
B_{k,0} &= -\frac{1}{k} \puz (A_{k-1,0} + D_{k-1,0}), \ \ k \geq 1\\
B_{k,\ell} &= (-1)^k \frac{(\ell-k)!}{k! \ell!} (\puz \puzd)^k B_{0,\ell-k}, \ \ 1 \leq k \leq \ell\\
B_{k, \ell} &= (-1)^{\ell+1} \frac{(k-1-\ell)!}{k! \ell!} \puz (\puzd \puz)^{\ell} (A_{k-1-\ell,0} + D_{k-1-\ell,0}), \ \ 1 \leq \ell < k
\end{align*}
and 
\begin{align*}
A_{0,\ell} + D_{0,\ell} &= \frac{1}{\ell} \puzd B_{0,\ell-1}, \ \ \ell \geq 1\\
A_{k,\ell} + D_{k,\ell} &= (-1)^{k} \frac{(\ell-1-k)!}{k! \ell!} \puzd (\puz \puzd)^k B_{0,\ell-1-k}, \ \ 1 \leq k < \ell\\
A_{k,\ell} + D_{k,\ell} &= (-1)^{\ell} \frac{(k-\ell)!}{k! \ell!} (\puzd \puz)^{\ell} (A_{k-\ell,0} + D_{k-\ell,0}), \ \ 1 \leq \ell \leq k.
\end{align*}
Hence these expressions allow us to determine $A_{k,\ell}$, $B_{k,\ell}$, $C_{k,\ell}$ and $D_{k,\ell}$, when $A_{0,\ell}$, $B_{0,\ell}$, $C_{k,0}$ and $D_{k,0}$ are given.

The double power series method is partitioned into three independent classes of solutions.
%%%%%%%%%%%%%%%%%%%%%%%%%%%%%%%%%%%%%%%%%%%%%%%%%%%%%%%%%%%%%%%%%%%%%%%%%%%%%%%%%%%%%%%%%%%%%%%%%%%%%%%%%%%%%%%%%%%%%%%%%%%%%%%%%%%%%%%%%%%%%%%%%%%%%%%%%%%%%%%%%%
\subsubsection{Class I}
This class was considered in subsection \ref{PSMet}. The functions $A$, $B$, $C$, $D$ then take the form (\ref{middle}).
%%%%%%%%%%%%%%%%%%%%%%%%%%%%%%%%%%%%%%%%%%%%%%%%%%%%%%%%%%%%%%%%%%%%%%%%%%%%%%%%%%%%%%%%%%%%%%%%%%%%%%%%%%%%%%%%%%%%%%%%%%%%%%%%%%%%%%%%%%%%%%%%%%%%%%%%%
\subsubsection{Class II}
In this case we put
\begin{equation}\label{first}
\begin{split}
A&=z_0^s \sum_{k=0}^\infty(z_0\overline z_0)^kA_k,\qquad\qquad B= z_0^{s+1} \sum_{k=0}^\infty (z_0\overline z_0)^kB_k,\\
C&=z_0^{s-1} \sum_{k=0}^\infty (z_0\overline z_0)^kC_k,\qquad\quad D= z_0^s \sum_{k=0}^\infty(z_0\overline z_0)^kD_k,
\end{split}
\end{equation}
with $s\in\mathbb N$ and $A_k$, $B_k$, $C_k$, $D_k$ $\mathbb{C}_{2n}$-valued continuously differentiable functions defined on $\mathbb R^{2n}$.
\begin{rem}
In case of $s>1$ we have that $f\vert_{z_0=0} = 0$, while in case of $s=1$ we obtain $f\vert_{z_0=0} = f_0^{\dagger} C_0$. However the restriction $f\vert_{z_0=0}$ does not determine the function $f$ uniquely.
\end{rem}
Plugging expressions (\ref{first}) in the systems of equations (\ref{whms1}) and (\ref{whms2}), we obtain the following recurrence relations
\begin{equation}\label{recf1}
A_{k} =\frac{1}{s+k}\partial_{\underline z}C_k, \ \ 
C_{k+1} =-\frac{1}{k+1}\partial_{\underline z^{\dagger}}A_k, \ \ k \geq 0
\end{equation}
and
\begin{equation}\label{recf2}
B_k =-\frac{1}{s+1+k} \partial_{\underline z}(A_k+D_k), \ \ 
A_{k+1}+D_{k+1} =\frac{1}{k+1}\partial_{\underline z^{\dagger}}B_k, \ \ k \geq 0.
\end{equation}
Next, (\ref{recf1}) leads to
\begin{align*}
A_k&= (-1)^k \frac{(s-1)!}{k! (s+k)!} \puz (\puzd \puz)^k C_0, \ \ k \geq 0\\
C_k&=  (-1)^k \frac{(s-1)!}{k! (s-1+k)!} (\puzd \puz)^k C_0, \ \ k \geq 1 
\end{align*}
while using (\ref{recf2}) we get
\begin{align*}
B_k&= (-1)^{k+1} \frac{s!}{k! (s+1+k)!} \puz (\puzd \puz)^k D_0, \ \ k \geq 0\\
D_k&= (-1)^k \frac{(s-1)!}{k! (s+k)!} \left( s (\puzd \puz)^k D_0 - \puz (\puzd \puz)^k C_0 \right), \ \ k \geq 1.
\end{align*}
Hence, the above equations allow us to determine all $A_k$, $B_k$, $C_k$ and $D_k$ given the starting values $C_0$ and $D_0$.

If we take for example as initial conditions the Gauss function, i.e. $C_0 = D_0 = e^{-\frac{\vert \uz \vert^2}{2}}$, the solutions can again be expressed in terms of the Hermitian Clifford-Hermite polynomials as before or, alternatively, in terms of the Laguerre polynomials:
\begin{align*}
A_k &= (-1)^k 4^{-k} \frac{(s-1)!}{k! (s+k)!} e^{-\frac{\vert \uz \vert^2}{2}} H_{2k+1}^{(2)}(\uz,\uzd)\\
&= -2^{-k-1} \frac{(s-1)!}{(s+k)!} e^{-\frac{\vert \uz \vert^2}{2}} \uzd  L_{k}^n \left( \frac{\vert \uz \vert^2}{2} \right), \ k \geq 0
\end{align*}
\begin{align*}
B_k &= (-1)^{k+1} 4^{-k}\frac{s!}{k! (s+1+k)!} e^{-\frac{\vert \uz \vert^2}{2}} H_{2k+1}^{(2)}(\uz,\uzd)\\
 &= 2^{-k-1} \frac{s!}{(s+1+k)!}  e^{-\frac{\vert \uz \vert^2}{2}} \uzd L_{k}^n \left( \frac{\vert \uz \vert^2}{2} \right), \ k \geq 0
\end{align*}
\begin{align*}
C_k&= (-1)^{k} 4^{1-k}\frac{(s-1)!}{k! (s+k-1)!} e^{-\frac{\vert \uz \vert^2}{2}} H_{2k}^{(4)}(\uz,\uzd)\\
&=2^{-k} \frac{(s-1)!}{k (s+k-1)!}  e^{-\frac{\vert \uz \vert^2}{2}} \left( (n-\beta) L_{k-1}^n \left( \frac{\vert \uz \vert^2}{2} \right) - \frac{1}{2} \uz\,\uzd  L_{k-1}^{n+1} \left( \frac{\vert \uz \vert^2}{2} \right) \right) , \ k \geq 1
\end{align*}
\begin{align*}
D_k&= (-1)^k 4^{-k} \frac{(s-1)!}{k! (s+k)!} e^{-\frac{\vert \uz \vert^2}{2}}  \left( 4s H_{2k}^{(4)}(\uz,\uzd) - H_{2k+1}^{(2)}(\uz,\uzd) \right)\\
&= 2^{-k-1} \frac{(s-1)!}{(s+k)!} e^{-\frac{\vert \uz \vert^2}{2}} \biggl\lbrace \frac{2s}{k} \left( (n-\beta)  L_{k-1}^n \left( \frac{\vert \uz \vert^2}{2} \right) - \frac{1}{2} \uz\,\uzd L_{k-1}^{n+1} \left( \frac{\vert \uz \vert^2}{2} \right)  \right)\\
&\qquad+\uzd L_{k}^{n} \left( \frac{\vert \uz \vert^2}{2} \right)  \biggr\rbrace, \ k \geq 1.
\end{align*}
%%%%%%%%%%%%%%%%%%%%%%%%%%%%%%%%%%%%%%%%%%%%%%%%%%%%%%%%%%%%%%%%%%%%%%%%%%%%%%%%%%%%%%%%%%%%%%%%%%%%%%%%%%%%%%%%%%%%%%%%%%%%%%%%%%%%%%%%%%%%%%%%%%%%%%%%%%%%%%%%%%%%%%%%%%%%%%%%%%%%%
\subsubsection{Class III}
Here we look for solutions of the systems (\ref{whms1}) and (\ref{whms2}) of the form
\begin{equation}\label{third}
\begin{split}
A&=\overline{z}_0^s \sum_{k=0}^\infty(z_0\overline z_0)^kA_k,\qquad\qquad B=\overline{z}_0^{s-1} \sum_{k=0}^\infty (z_0\overline z_0)^kB_k,\\
C&=\overline{z}_0^{s+1} \sum_{k=0}^\infty(z_0\overline z_0)^kC_k,\qquad\quad D=\overline{z}_0^s \sum_{k=0}^\infty(z_0\overline z_0)^kD_k,
\end{split}
\end{equation}
with again $s\in\mathbb N$ and $A_k$, $B_k$, $C_k$, $D_k$ $\mathbb{C}_{2n}$-valued continuously differentiable functions defined on $\mathbb R^{2n}$.
\begin{rem}
We have that
\begin{equation*}
f\vert_{z_0=0} = 
\begin{cases}
0 & \text{if $s>1$},\\
f_0 B_0 & \text{if $s=1$}.
\end{cases}
\end{equation*}
However, again this restriction does not determine the function $f$ uniquely.
\end{rem}
Substituting the expressions (\ref{third}) in the systems (\ref{whms1}) and (\ref{whms2}), we arrive at the following recurrence relations:
\begin{equation}\label{rect1}
A_{k+1} =\frac{1}{k+1} \partial_{\underline z}C_k, \ \ 
C_{k} =-\frac{1}{s+1+k} \partial_{\underline z^{\dagger}} A_k, \ \ k \geq 0
\end{equation}
and
\begin{equation}\label{rect2}
B_{k+1} =-\frac{1}{k+1}\partial_{\underline z}(A_k + D_k), \ \ 
A_{k}+D_{k} =\frac{1}{s+k}\partial_{\underline z^{\dagger}}B_k, \ \ k \geq 0.
\end{equation} 
From (\ref{rect1}) we obtain expressions for $A_k$ and $C_k$ in terms of $A_0$:
\begin{align*}
A_k&= (-1)^k \frac{s!}{k! (s+k)!} (\puz \puzd)^k A_0, \ \ k \geq 1\\
C_k&=  (-1)^{k+1} \frac{s!}{k! (s+k+1)!} \puzd (\puz \puzd)^k A_0, \ \ k \geq 0 
\end{align*}
while (\ref{rect2}) yields expressions for $B_k$ and $D_k$ in terms of the starting values $A_0$ and $B_0$:
\begin{align*}
B_k&= (-1)^{k} \frac{(s-1)!}{k! (s+k-1)!} (\puz \puzd)^k B_0, \ \ k \geq 1\\
D_k&= (-1)^k \frac{(s-1)!}{k! (s+k)!} \left( \puzd (\puz \puzd)^k B_0 - s (\puz \puzd)^k A_0 \right), \ \ k \geq 0.
\end{align*}
Let us again illustrate the above with the Gauss function as initial conditions, i.e.  $A_0 = B_0 = e^{-\frac{\vert \uz \vert^2}{2}}$. Like before, the solutions can be expressed in terms of the Hermitian Clifford-Hermite polynomials, or, alternatively, the Laguerre polynomials:
\begin{align*}
A_k &= (-1)^k 4^{1-k} \frac{s!}{k! (s+k)!} e^{-\frac{\vert \uz \vert^2}{2}} H_{2k}^{(3)}(\uz,\uzd)\\
& = 2^{-k} \frac{s!}{k (s+k)!} e^{-\frac{\vert \uz \vert^2}{2}} \left( \beta L_{k-1}^n \left( \frac{\vert \uz \vert^2}{2} \right) - \frac{1}{2} \uzd \uz  L_{k-1}^{n+1} \left( \frac{\vert \uz \vert^2}{2} \right)  \right), \ k \geq 1
\end{align*}
\begin{align*}
B_k &= (-1)^{k} 4^{1-k}\frac{(s-1)!}{k! (s+k-1)!} e^{-\frac{\vert \uz \vert^2}{2}} H_{2k}^{(3)}(\uz,\uzd)\\
& = 2^{-k} \frac{(s-1)!}{k (s+k-1)!}  e^{-\frac{\vert \uz \vert^2}{2}}  \left( \beta L_{k-1}^n \left( \frac{\vert \uz \vert^2}{2} \right) - \frac{1}{2} \uzd \uz L_{k-1}^{n+1} \left( \frac{\vert \uz \vert^2}{2} \right) \right), \ k \geq 1
\end{align*}
\begin{align*}
C_k &= (-1)^{k+1} 4^{-k}\frac{s!}{k! (s+k+1)!} e^{-\frac{\vert \uz \vert^2}{2}} H_{2k+1}^{(1)}(\uz,\uzd)\\
&= 2^{-k-1} \frac{s!}{(s+k+1)!}  e^{-\frac{\vert \uz \vert^2}{2}} \uz L_{k}^n \left( \frac{\vert \uz \vert^2}{2} \right) , \ k \geq 0
\end{align*}
\begin{align*}
D_k &= (-1)^k 4^{-k} \frac{(s-1)!}{k! (s+k)!} e^{-\frac{\vert \uz \vert^2}{2}}  \left( H_{2k+1}^{(1)}(\uz,\uzd) - 4s H_{2k}^{(3)}(\uz,\uzd) \right)\\
&= - 2^{-k-1} \frac{(s-1)!}{(s+k)!} e^{-\frac{\vert \uz \vert^2}{2}} \biggl\lbrace \uz  L_{k}^n \left( \frac{\vert \uz \vert^2}{2} \right) + \frac{2s}{k} \biggl( \beta L_{k-1}^{n} \left( \frac{\vert \uz \vert^2}{2} \right)\\
&\qquad-\frac{1}{2} \uzd \uz  L_{k-1}^{n+1}  \left( \frac{\vert \uz \vert^2}{2} \right) \biggr) \biggr\rbrace, \ k \geq 1\\
D_0 &= - e^{-\frac{\vert \uz \vert^2}{2}} \left( \frac{1}{2s} \uz + 1 \right).
\end{align*}
%%%%%%%%%%%%%%%%%%%%%%%%%%%%%%%%%%%%%%%%%%%%%%%%%%%%%%%%%%%%%%%%%%%%%%%%%%%%%%%%%%%%%%%%%%%%%%%%%%%%%%%%%%%%%%%%%%%%%%%%%%%%%%%%%%%%%%%%%%%
\section{Axial-type solutions}
\subsection{The Hermitian Vekua system}\label{Vekua}
Here, we will look for special solutions of the systems (\ref{whms1}) and (\ref{whms2}), which we assume to be of the form
\begin{align*}
A&=a_1+\underline z^{\dagger}\underline za_2,&
B&=z_0\underline z^{\dagger}b,\\
C&=\overline z_0\underline zc,&
D&=d_1+\underline z^{\dagger}\underline zd_2,
\end{align*}
where $a_1$, $a_2$, $b$, $c$, $d_1$, $d_2$ are continuously differentiable functions 
depending on the two variables 
\[(\nu_0,\nu)=(\vert z_0\vert^2,\vert \uz \vert^2),\]
and taking values in the real algebra generated by $\beta=\sum_{j=1}^nf_j^\dagger f_j$, which, as we know, is the quotient of the free ring generated by $\beta$ with its 2-sided ideal generated by the polynomial $\beta(\beta-1)\dots(\beta-n)$. One may also look for spinor valued solutions in $\Lambda_n^{\ell}I$ thus allowing $\beta$ to be replaced by its eigenvalue $\ell$.

It is easily seen that
\begin{align*}
\partial_{z_0}A&=\overline z_0\left(\partial_{\nu_0}a_1+\underline z^{\dagger}\underline z\partial_{\nu_0}a_2\right)\\
\partial_{\underline z}C&=\overline z_0\left((\partial_{\underline z}\underline z)c+\sum_{j=1}^nf_j^\dagger\underline z(\partial_{z_j}c)\right)=\overline z_0\left(\beta c+\underline z^{\dagger}\underline z\partial_{\nu}c\right)\\
\partial_{\overline z_0}C&=\underline z(c+\nu_0\partial_{\nu_0}c).
\end{align*}
Using the identity 
\[\underline z\beta-\beta\underline z=\underline z,\]
we also obtain that
\begin{align*}
\partial_{\underline z^{\dagger}}A&=\underline z\partial_{\nu}a_1+\sum_{j=1}^nf_j\left(f_j^\dagger\underline za_2+z_j\underline z^{\dagger}\underline z\partial_{\nu}a_2\right)=\underline z\partial_{\nu}a_1+(n-\beta)\underline za_2+\underline z\,\underline z^\dagger\underline z\partial_{\nu}a_2\\
&=\underline z\left(\partial_{\nu}a_1+(n+1-\beta)a_2+\nu\partial_{\nu}a_2\right).
\end{align*}
Therefore, system (\ref{whms1}) takes the form
\begin{equation}\label{whms1axial}
\left\{\begin{aligned}
\partial_{\nu_0}a_1&=\beta c\\
\partial_{\nu_0}a_2-\partial_{\nu}c&=0\\
\partial_{\nu}a_1+\nu\partial_{\nu}a_2+\nu_0\partial_{\nu_0}c&=(\beta-n-1)a_2-c.
\end{aligned}\right.
\end{equation}
Making similar calculations, we may check that system (\ref{whms2}) can be rewritten as
\begin{equation}\label{whms2axial}
\left\{\begin{aligned}
\partial_{\nu}(a_1+d_1)+\nu_0\partial_{\nu_0}b&=\beta(a_2+d_2)-b\\
\partial_{\nu_0}(a_1+d_1)-\nu\partial_\nu b&=(n-\beta) b\\
\partial_{\nu_0}(a_2+d_2)+\partial_{\nu}b&=0.
\end{aligned}\right.
\end{equation}
\begin{rem}
Note that systems (\ref{whms1axial}) and (\ref{whms2axial}) are the Hermitian equivalent of the Vekua system considered in \cite{LB,S1,S2} that describes axially symmetric monogenic functions in the orthogonal setting.
\end{rem}
In order to solve the Vekua systems (\ref{whms1axial}) and (\ref{whms2axial}), we shall make use of the power series method by writing
\begin{align*}
a_j(\nu_0,\nu)&=\sum_{k=0}^\infty\nu_0^ka_{k,j}(\nu),&
b(\nu_0,\nu)&=\sum_{k=0}^\infty\nu_0^kb_{k}(\nu),\\
c(\nu_0,\nu)&=\sum_{k=0}^\infty\nu_0^kc_{k}(\nu),&
d_j(\nu_0,\nu)&=\sum_{k=0}^\infty\nu_0^kd_{k,j}(\nu),
\end{align*}
with $j=1,2$. The above systems may now be rewritten in the form
\begin{align}
a_{k,1}&=\frac{\beta}{k}c_{k-1},&k\ge1\label{eq1}\\
a_{k,2}&=\frac{c_{k-1}^\prime}{k},&k\ge1\label{eq2}\\
(k+1)c_k+a_{k,1}^\prime+\nu a_{k,2}^\prime&=(\beta-n-1)a_{k,2},&k\ge0\label{eq3}
\end{align}
\begin{align}
(k+1)b_k+a_{k,1}^\prime+d_{k,1}^\prime&=\beta(a_{k,2}+d_{k,2}),&k\ge0\label{eq4}\\
a_{k,1}+d_{k,1}&=\frac{1}{k}\left((n-\beta)b_{k-1}+\nu b_{k-1}^\prime\right),&k\ge1\label{eq5}\\
a_{k,2}+d_{k,2}&=-\frac{b_{k-1}^\prime}{k},&k\ge1.\label{eq6}
\end{align}
Substituting (\ref{eq1}), (\ref{eq2}) in (\ref{eq3}), we get
\begin{equation}\label{rraxial1}
\begin{split}
c_k&=-\frac{1}{k(k+1)}\left((n+1)c_{k-1}^\prime+\nu c_{k-1}^{\prime\prime}\right),\quad k\ge1\\
c_0&=(\beta-n-1)a_{0,2}-a_{0,1}^\prime-\nu a_{0,2}^\prime.
\end{split}
\end{equation}
In a similar way, substituting (\ref{eq5}), (\ref{eq6}) in (\ref{eq4}), we obtain that
\begin{equation}\label{rraxial2}
\begin{split}
b_k&=-\frac{1}{k(k+1)}\left((n+1)b_{k-1}^\prime+\nu b_{k-1}^{\prime\prime}\right),\quad k\ge1\\
b_0&=\beta(a_{0,2}+d_{0,2})-a_{0,1}^\prime-d_{0,1}^\prime.
\end{split}
\end{equation}
Finally, from (\ref{eq5}), (\ref{eq6}) and using (\ref{eq1}), (\ref{eq2}) we get
\begin{align}
d_{k,1}&=\frac{1}{k}\left((n-\beta)b_{k-1}+\nu b_{k-1}^\prime-\beta c_{k-1}\right),&k\ge1\label{rraxial3}\\
d_{k,2}&=-\frac{1}{k}(b_{k-1}^\prime+c_{k-1}^\prime),&k\ge1.\label{rraxial4}
\end{align}
Hence the explicit solutions of the Vekua systems (\ref{whms1axial}), (\ref{whms2axial}) may be obtained via the recurrence relations (\ref{rraxial1}), (\ref{rraxial2}), (\ref{eq1}), (\ref{eq2}), (\ref{rraxial3}) and (\ref{rraxial4}) with the initial conditions 
\[a_j(0,\nu),\;d_j(0,\nu),\quad j=1,2.\]
%%%%%%%%%%%%%%%%%%%%%%%%%%%%%%%%%%%%%%%%%%%%%%%%%%%%%%%%%%%%%%%%%%%%%%%
\subsection{Hermitian generalized powers}
Let us illustrate the above ideas with an example that relates to the definition of generalized powers in the Hermitian setting. For the orthogonal setting we refer the reader to \cite{S2}). Suppose that
\[a_{0,1}(\nu)=\alpha_1\nu^s,\quad a_{0,2}(\nu)=\alpha_2\nu^{s-1},\]
\[d_{0,1}(\nu)=\delta_1\nu^s,\quad d_{0,2}(\nu)=\delta_2\nu^{s-1},\]
with $s,\alpha_j,\delta_j\in\mathbb R$, $j=1,2$. Using the previous recurrence relations we easily obtain
\[b_k(\nu)=\alpha(k)\nu^{s-k-1},\;\;k\ge0,\]
\[c_k(\nu)=\delta(k)\nu^{s-k-1},\;\;k\ge0,\]
where
\[\alpha(k)=\frac{(-1)^k\prod_{\ell=1}^k(s-\ell)\prod_{\ell=1}^k(n+s-\ell) \left(\beta(\alpha_2+\delta_2)-(\alpha_1+\delta_1)s\right)}{(k+1)(k!)^2},\]
\[\delta(k)=\frac{(-1)^k\prod_{\ell=1}^k(s-\ell)\prod_{\ell=1}^k(n+s-\ell)\left((\beta-n-s)\alpha_2-\alpha_1s\right)}{(k+1)(k!)^2}.\]
Hence
\[a_{k,1}(\nu)=\frac{\beta\delta(k-1)}{k}\nu^{s-k},\]
\[a_{k,2}(\nu)=\frac{(s-k)\delta(k-1)}{k}\nu^{s-k-1},\]
\[d_{k,1}(\nu)=\frac{(n+s-k-\beta)\alpha(k-1)-\beta\delta(k-1)}{k}\nu^{s-k},\]
\[d_{k,2}(\nu)=-\frac{(s-k)}{k}\left(\alpha(k-1)+\delta(k-1)\right)\nu^{s-k-1},\]
with $k\ge1$.
\begin{rem}
It is worth pointing out that if $s\in\mathbb N$ then $\alpha(k)=\delta(k)=0$ for $k\ge s$. Therefore, for the case $s\in\mathbb N$ the solutions $A$, $B$, $C$, $D$ are homogeneous polynomials of degree $2s$ in $\mathbb R^{2n+2}$.
\end{rem}
%%%%%%%%%%%%%%%%%%%%%%%%%%%%%%%%%%%%%%%%%%%%%%%%%%%%%%%%%%%%%%%%%%%%%%
\subsection{Two more cases}
We may consider two more cases of axial-type solutions of the systems (\ref{whms1}) and (\ref{whms2}). Indeed, if we assume that
\begin{align*}
A&=z_0^s(a_1+\underline z^{\dagger}\underline za_2),&
B&=z_0^{s+1}\underline z^{\dagger}b,\\
C&=z_0^{s-1}\underline zc,&
D&=z_0^s(d_1+\underline z^{\dagger}\underline zd_2),
\end{align*}
with $s\in\mathbb N$, then we obtain the following systems of equations
\begin{equation*}
\left\{\begin{aligned}
\nu_0\partial_{\nu_0}a_1&=-sa_1+\beta c\\
\nu_0\partial_{\nu_0}a_2-\partial_{\nu}c&=-sa_2\\
\partial_{\nu}a_1+\nu\partial_{\nu}a_2+\partial_{\nu_0}c&=(\beta-n-1)a_2
\end{aligned}\right.
\end{equation*}
\begin{equation*}
\left\{\begin{aligned}
\partial_{\nu}(a_1+d_1)+\nu_0\partial_{\nu_0}b&=\beta(a_2+d_2)-(s+1)b\\
\partial_{\nu_0}(a_1+d_1)-\nu\partial_\nu b&=(n-\beta) b\\
\partial_{\nu_0}(a_2+d_2)+\partial_{\nu}b&=0
\end{aligned}\right.
\end{equation*}
which can be solved using the power series method as we did in subsection \ref{Vekua}. The explicit solutions are given as follows
\[c_k=-\frac{1}{k(k+s-1)}\left((n+1)c_{k-1}^\prime+\nu c_{k-1}^{\prime\prime}\right),\quad k\ge1\]
\begin{align*}
a_{k,1}&=\frac{\beta}{k+s}\,c_{k},&k\ge0\\
a_{k,2}&=\frac{c_{k}^\prime}{k+s},&k\ge0
\end{align*}
\begin{equation*}
\begin{split}
b_k&=-\frac{1}{k(k+s+1)}\left((n+1)b_{k-1}^\prime+\nu b_{k-1}^{\prime\prime}\right),\quad k\ge1\\
b_0&=\frac{1}{s+1}(\beta d_{0,2}-d_{0,1}^\prime)
\end{split}
\end{equation*}
\begin{align*}
d_{k,1}&=\frac{1}{k}\left((n-\beta)b_{k-1}+\nu b_{k-1}^\prime\right)-\frac{\beta}{k+s}\,c_{k},&k\ge1\\
d_{k,2}&=-\frac{b_{k-1}^\prime}{k}-\frac{c_{k}^\prime}{k+s},&k\ge1
\end{align*}
with the initial conditions 
\[c(0,\nu),\;d_j(0,\nu),\quad j=1,2.\]
Finally, let us assume that
\begin{align*}
A&=\overline z_0^s(a_1+\underline z^{\dagger}\underline za_2),&
B&=\overline z_0^{s-1}\underline z^{\dagger}b,\\
C&=\overline z_0^{s+1}\underline zc,&
D&=\overline z_0^s(d_1+\underline z^{\dagger}\underline zd_2),
\end{align*}
with $s\in\mathbb N$. This yields the systems
\begin{equation*}
\left\{\begin{aligned}
\partial_{\nu_0}a_1&=\beta c\\
\partial_{\nu_0}a_2-\partial_{\nu}c&=0\\
\partial_{\nu}a_1+\nu\partial_{\nu}a_2+\nu_0\partial_{\nu_0}c&=(\beta-n-1)a_2-(s+1)c
\end{aligned}\right.
\end{equation*}
\begin{equation*}
\left\{\begin{aligned}
\partial_{\nu}(a_1+d_1)+\partial_{\nu_0}b&=\beta(a_2+d_2)\\
\nu_0\partial_{\nu_0}(a_1+d_1)-\nu\partial_\nu b&=-s(a_1+d_1)+(n-\beta) b\\
\nu_0\partial_{\nu_0}(a_2+d_2)+\partial_{\nu}b&=-s(a_2+d_2),
\end{aligned}\right.
\end{equation*}
whose solutions may be obtained via the recurrence relations  
\begin{equation*}
\begin{split}
c_k&=-\frac{1}{k(k+s+1)}\left((n+1)c_{k-1}^\prime+\nu c_{k-1}^{\prime\prime}\right),\quad k\ge1\\
c_0&=\frac{1}{s+1}\left((\beta-n-1)a_{0,2}-a_{0,1}^\prime-\nu a_{0,2}^\prime\right)
\end{split}
\end{equation*}
\begin{align*}
a_{k,1}&=\frac{\beta}{k}\,c_{k-1},&k\ge1\\
a_{k,2}&=\frac{c_{k-1}^\prime}{k},&k\ge1
\end{align*}
\begin{equation*}
b_k=-\frac{1}{k(k+s-1)}\left((n+1)b_{k-1}^\prime+\nu b_{k-1}^{\prime\prime}\right),\quad k\ge1
\end{equation*}
\begin{align*}
d_{k,1}&=\frac{1}{k+s}\left((n-\beta)b_{k}+\nu b_{k}^\prime\right)-a_{k,1},&k\ge0\\
d_{k,2}&=-\frac{b_{k}^\prime}{k+s}-a_{k,2},&k\ge0
\end{align*}
with the initial conditions 
\[a_j(0,\nu),\;b(0,\nu),\quad j=1,2.\]
%%%%%%%%%%%%%%%%%%%%%%%%%%%%%%%%%%%%%%%%%%%%%%%%%%%%%%%%%%%%%%%%%%%%%%%%%%%%%%%%%%%%%%%%%%%%%%%%%%%%%%%%%%%%%%%%%%%%%%%%%%%%%%%%%%%%%%
\section{Exponential-type solutions and Hermitian Bessel functions}
Here we shall seek special solutions of the systems (\ref{whms1}) and (\ref{whms2}) of the form
\begin{align*}
A&=e^{\lambda z_0+\mu\overline z_0}(a_1+\underline z^{\dagger}\underline za_2),&
B&=e^{\lambda z_0+\mu\overline z_0}\underline z^{\dagger}b,\\
C&=e^{\lambda z_0+\mu\overline z_0}\underline zc,&
D&=e^{\lambda z_0+\mu\overline z_0}(d_1+\underline z^{\dagger}\underline zd_2),
\end{align*}
$\lambda,\mu\in\mathbb R\setminus\{0\}$ and where $a_1$, $a_2$, $b$, $c$, $d_1$, $d_2$ are continuously differentiable functions depending on the variable $\nu=\vert\uz\vert^2$ taking values in the real algebra generated by $\beta$.

In what follows $J_\alpha$, $I_\alpha$ stand, respectively, for the Bessel function of the first kind and the modified Bessel function (see e.g. \cite{Hoch}), given by
\[J_\alpha(t)=\sum_{k=0}^\infty\frac{(-1)^k}{k!\,\Gamma(k+\alpha+1)} {\left(\frac{t}{2}\right)}^{2k+\alpha},\]
\[I_\alpha(t)=\sum_{k=0}^\infty \frac{1}{k!\,\Gamma(k+\alpha+1)}\left(\frac{t}{2}\right)^{2k+\alpha}.\]
Proceeding in the same spirit as in the previous section we can check that systems (\ref{whms1}) and (\ref{whms2}) now take the form 

\begin{equation*}\label{whms1exp}
\left\{\begin{aligned}
a_1&=\frac{\beta}{\lambda}c\\
a_2&=\frac{c^\prime}{\lambda}\\
a_1^\prime+\nu a_2^\prime+\mu c&=(\beta-n-1)a_2
\end{aligned}\right.
\end{equation*}
\begin{equation*}\label{whms2exp}
\left\{\begin{aligned}
a_1^\prime+d_1^\prime+\lambda b&=\beta(a_2+d_2)\\
a_1+d_1&=\frac{1}{\mu}\left((n-\beta)b+\nu b^\prime\right)\\
a_2+d_2&=-\frac{b^\prime}{\mu}.
\end{aligned}\right.
\end{equation*}
From these systems of equations we can easily deduce that $b$ and $c$ satisfy the same ordinary differential equation, namely
\begin{align*}
\nu b^{\prime\prime}+(n+1)b^\prime+\lambda\mu b&=0\\
\nu c^{\prime\prime}+(n+1)c^\prime+\lambda\mu c&=0.
\end{align*}
It has solutions of the form
\[b(\nu)=\left\{\begin{array}{ll}\alpha_1\nu^{-\frac{n}{2}}J_n\left(2\sqrt{\lambda\mu\nu}\right)&\text{if}\;\;\lambda\mu>0,\\
\alpha_1\nu^{-\frac{n}{2}}I_n\left(2\sqrt{-\lambda\mu\nu}\right)&\text{if}\;\;\lambda\mu<0,\end{array}\right.\]
\[c(\nu)=\left\{\begin{array}{ll}\alpha_2\nu^{-\frac{n}{2}}J_n\left(2\sqrt{\lambda\mu\nu}\right)&\text{if}\;\;\lambda\mu>0,\\
\alpha_2\nu^{-\frac{n}{2}}I_n\left(2\sqrt{-\lambda\mu\nu}\right)&\text{if}\;\;\lambda\mu<0,\end{array}\right.\]
where $\alpha_1$, $\alpha_2$ are arbitrary real constants. The other functions can be easily obtained from the above equations. Indeed,
\[a_1(\nu)=\frac{\beta}{\lambda}c(\nu),\]
\[a_2(\nu)=\left\{\begin{array}{ll}\displaystyle{-\alpha_2\frac{\sqrt{\lambda\mu}}{\lambda}\,\nu^{-\frac{n+1}{2}}J_{n+1}\left(2\sqrt{\lambda\mu\nu}\right)}&\text{if}\;\;\lambda\mu>0,\\
\displaystyle{\alpha_2\frac{\sqrt{-\lambda\mu}}{\lambda}\,\nu^{-\frac{n+1}{2}}I_{n+1}\left(2\sqrt{-\lambda\mu\nu}\right)}&\text{if}\;\;\lambda\mu<0,\end{array}\right.\]
\[d_2(\nu)=\left\{\begin{array}{ll}\displaystyle{\sqrt{\lambda\mu}\left(\frac{\alpha_1}{\mu}+\frac{\alpha_2}{\lambda}\right)\nu^{-\frac{n+1}{2}}J_{n+1}\left(2\sqrt{\lambda\mu\nu}\right)}&\text{if}\;\;\lambda\mu>0,\\
\displaystyle{-\sqrt{-\lambda\mu}\left(\frac{\alpha_1}{\mu}+\frac{\alpha_2}{\lambda}\right)\nu^{-\frac{n+1}{2}}I_{n+1}\left(2\sqrt{-\lambda\mu\nu}\right)}&\text{if}\;\;\lambda\mu<0,\end{array}\right.\]
\[d_1(\nu)=\frac{1}{\mu}(n-\beta)b(\nu)-\frac{\beta}{\lambda}c(\nu)-\nu\big(a_2(\nu)+d_2(\nu)\big).\]
%%%%%%%%%%%%%%%%%%%%%%%%%%%%%%%%%%%%%%%%%%%%%%%%%%%%%%%%%%%%%%%%%%%%%%%%%%%%%%%%%%%%%%%%%%%%%%%%%%%%%%%%%%%%%%%%%%%%%%%%%%%%%%%%%%%%%%%%%%%%%%%%%%%%
\subsection*{Acknowledgments}
D. Pe\~na Pe\~na acknowledges the support of a Postdoctoral Fellowship from \lq\lq Special Research Fund" (BOF) of Ghent University.
%%%%%%%%%%%%%%%%%%%%%%%%%%%%%%%%%%%%%%%%%%%%%%%%%%%%%%%%%%%%%%%%%%%%%%%

\end{document}